\documentclass[fleqn,11pt,twoside]{article}

\usepackage{amsthm,amsthm,amssymb, color, xcolor,epsfig, graphics, subfigure}

\usepackage{amsmath, graphicx, latexsym, lscape}

\usepackage{cite}
\usepackage{float}

\usepackage[shortlabels]{enumitem}

\makeatletter
\newcommand{\copyrightnote}[2]{{\renewcommand{\thefootnote}{}
 \footnotetext{\small\it
\begin{flushleft}
 \copyright \ #1   #2  
\end{flushleft}}}}

\newcommand{\Name}[1]{\begin{flushleft}
                       \LARGE \bf #1
                       \end{flushleft}\vspace{-3mm}}

\newcommand{\Author}[1]{\begin{flushleft}
                       \it #1 \end{flushleft}}

\newcommand{\Address}[1]{\begin{flushleft}
                       \it #1 \end{flushleft}}

\newcommand{\Date}[1]{\begin{flushleft}
                      \small  \it #1 \end{flushleft}}

\newcommand{\evenhead}{Author \ name}
\newcommand{\oddhead}{Article \ name}

\renewcommand{\@evenhead}{
\hspace*{-3pt}\raisebox{-15pt}[\headheight][0pt]{\vbox{\hbox to \textwidth
{\thepage \hfil \evenhead}\vskip4pt \hrule}}}
\renewcommand{\@oddhead}{
\hspace*{-3pt}\raisebox{-15pt}[\headheight][0pt]{\vbox{\hbox to \textwidth
{\oddhead \hfil \thepage}\vskip4pt\hrule}}}
\renewcommand{\@evenfoot}{}
\renewcommand{\@oddfoot}{}

\long\def\@makecaption#1#2{%
  \vskip\abovecaptionskip
  \sbox\@tempboxa{\small \textbf{#1.}\ \ #2}%
  \ifdim \wd\@tempboxa >\hsize
    {\small \textbf{#1.}\ \ #2}\par
  \else
    \global \@minipagefalse
    \hb@xt@\hsize{\hfil\box\@tempboxa\hfil}%
  \fi
  \vskip\belowcaptionskip}

\newcommand{\JNMPnumberwithin}[3][\arabic]{%
  \@ifundefined{c@#2}{\@nocounterr{#2}}{%
    \@ifundefined{c@#3}{\@nocnterr{#3}}{%
      \@addtoreset{#2}{#3}%
      \@xp\xdef\csname the#2\endcsname{%
        \@xp\@nx\csname the#3\endcsname .\@nx#1{#2}}}}%
}

\newcommand{\resetfootnoterule} {
  \renewcommand\footnoterule{%
  \kern-3\p@
  \hrule\@width.4\columnwidth
  \kern2.6\p@}
}

\renewcommand{\footnoterule}{}

\makeatother

\theoremstyle{definition}

\newcommand{\itr}[2]{\overset{(#1)}{#2}}

\begin{document}

\renewcommand{\evenhead}{ {\LARGE\textcolor{blue!10!black!40!green}{{\sf \ \ \ ]ocnmp[}}}\strut\hfill E I Kaptsov and V A Dorodnitsyn}
\renewcommand{\oddhead}{ {\LARGE\textcolor{blue!10!black!40!green}{{\sf ]ocnmp[}}}\ \ \ \ \  Invariant conservative finite-difference schemes}

\thispagestyle{empty}
\newcommand{\FistPageHead}[3]{
\begin{flushleft}
\raisebox{8mm}[0pt][0pt]
{\footnotesize \sf
\parbox{150mm}{{Open Communications in Nonlinear Mathematical Physics}\ \ \ \ {\LARGE\textcolor{blue!10!black!40!green}{]ocnmp[}}
\quad Special Issue 1, 2024\ \  pp
#2\hfill {\sc #3}}}\vspace{-13mm}
\end{flushleft}}

\FistPageHead{1}{\pageref{firstpage}--\pageref{lastpage}}{ \ \ }

\strut\hfill

\strut\hfill

\copyrightnote{The author(s). Distributed under a Creative Commons Attribution 4.0 International License}

\begin{center}
{  {\bf This article is part of an OCNMP Special Issue\\ 
\smallskip
in Memory of Professor Decio Levi}}
\end{center}

\smallskip

\Name{Invariant Conservative Finite-Difference Schemes For the 
One-Dimensional Shallow Water Magneto\-hydrodynamics Equations 
in~Lagrangian Coordinates}

\Author{E.~I. Kaptsov$^{\,1}$ and V.~A. Dorodnitsyn$^{\,2}$}

\Address{$^{1}$ School of Mathematics, Institute of Science, 
 Suranaree University of Technology, 30000, Thailand\\[2mm]
$^{2}$ Keldysh Institute of Applied Mathematics, 
 Russian Academy of Science, Miusskaya Pl.~4, Moscow, 125047, Russia}

\Date{Received April 26, 2023 ; Accepted September 14, 2023}

\setcounter{equation}{0}


\begin{abstract}
Finite-difference schemes for the one-dimensional shallow water
equations in the presence of a magnetic field for various bottom
topographies are constructed. Based on the results of the group
classification recently carried out by the authors,
finite-difference analogues of the conservation laws of the original
differential model are obtained. Some typical problems are
considered numerically, for which a comparison is made between the
cases of a magnetic field presence  and when it is absent (the
standard shallow water model). The invariance of difference schemes
in Lagrangian coordinates and the energy preservation on the
obtained numerical solutions are also discussed.
\end{abstract}


\label{firstpage}

\section{Introduction}

The consideration of the shallow water equations in the presence of a magnetic field (SMHD) is a relatively new area of magnetohydrodynamics. One of the first models to describe SMHD was presented in~\cite{bk:SWMHD_Gilman2000} (the Gilman model), and since then, many papers by various authors have been devoted to studying this model, including its stability~\cite{bk:SWMHD_Trakhinin2020}, numerical simulation~\cite{bk:SWMHD_Touma2010,bk:SWMHD_Winters2016,Bouchut2016,bk:SWMHD_Salem2017,bk:SWMHD_AhmedZia2019,bk:SWMHD_JunmingHuazhong2021}, and conservation laws \cite{bk:KapMelDor_PhA2022,bk:SWMHD_DeSterck2001} (see, for example, a more detailed review in~\cite{bk:KapMelDor_PhA2022}).

Initially, the SMHD equations were used to model the behavior of the
solar tachocline, which is a thin layer at the base of the Sun's
convection layer. Toroidal magnetic fields are usually considered in
this case, which can be assumed to act approximately in the
tachocline plane, so that the model~\cite{bk:SWMHD_Gilman2000}
describes two-dimensional flows. In some papers
(e.g.,~\cite{bk:SWMHD_SHIUE2013215,bk:KapMelDor_PhA2022}),
one-dimensional SMHD equations are considered.
In~\cite{bk:KapMelDor_PhA2022}, one-dimensional SMHD equations with
arbitrary bottom topography were studied in Lagrangian coordinates.
In this case, most of the equations of the SMHD system are
integrated, and only one equation remains unintegrated. For this
equation in~\cite{bk:KapMelDor_PhA2022}, a group classification was
carried out according to the function describing the topography of
the bottom. Conservation laws and invariant solutions in Lagrangian
and Eulerian coordinates were obtained. It was found that in
Lagrangian coordinates, the one-dimensional SMHD equations differ
from the standard one-dimensional shallow water
equations~\cite{bk:AndrKapPukhRod[1998],bk:KaptsovMeleshko_1D_classf[2018]}
only by a linear term. Thus, based on the known results for the
shallow water equations in Lagrangian coordinates, the
conservation laws for~SMHD was constructed, which can be done both
for differential equations and for their discretizations.

It is known~\cite{bk:Noether1918,bk:Ibragimov2[2007]} that the
presence of conservation laws for equations is closely related to
the symmetries of the equations. The inheritance of symmetries and
conservation laws is also important in the discretization of
equations and numerical
modeling~\cite{Dor_1,[LW-4],[LW-5],bk:LW2005,bk:Dorodnitsyn[2011]}.
The study of symmetries and conservation laws of discrete models,
including finite-difference equations, is the subject of a
significant number of works, such
as~\cite{levi_olver_thomova_winternitz_2011,Dor_1,[LW-4],[LW-5],[LW-6],bk:LW2005,bk:Dorodnitsyn[2011]}.
We would like to emphasize significant contributions of Professor Decio Levi to the study of 
discrete dynamical system symmetries~\cite{[LW-5]} and 
nonlinear differential difference equations as Bäcklund transformations, elucidating their connection to continuous evolution equations~\cite{DLevi_BacklundTr_1981}. 
D.~Levi has also stated conditions for the existence of higher symmetries~\cite{DLevi_HigherSyms_1997} and contact transformations~\cite{[LW-11],[LW-12]} for discrete equations.
Comprehensive research in the area of symmetries and integrability of discrete equations has been reflected in a recent book~\cite{bk:LeviWintYamilovBook2023} by D.~Levi, P.~Winternitz and I.~R.~Yamilov.

As part of the development of methods of group analysis of discrete equations, 
in~\cite{bk:DorodnitsynKozlovWint[2003],bk:Dorod_Hamilt[2011]},
finite-difference analogues of the~Lagrangian and~Hamiltonian
formalisms were proposed, which simplify the construction of 
invariant (symmetry-preserving) schemes and the derivation of their
conservation laws. If no invariant Lagrangian or Hamiltonian exists, 
the alternative method was introduced in~\cite{bk:DorodKozlovWintKaptsov[2015]}. 
Practically, to construct invariant schemes for
partial differential equations it is often more convenient to utilize 
the so-called direct
method~\cite{DorKap_2020,bk:ChevDorKap2020,bk:Bluman1997}, which was
recently used to obtain invariant conservative schemes for various
equations of hydrodynamic
type~\cite{DorKap_2020,bk:DorKapSW_JMP_2020,bk:KaptDorMel_ModSW2022,bk:DorKapMelGN2020,bk:DorKapMDPI2022,bk:KapDorMel_MHDCylSchemes2023,bk:KaptsovMHDCylSchemes_Conf2022}.

In the recent papers~\cite{DorKap_2020,bk:DorKapSW_JMP_2020} the authors constructed
finite-difference schemes for the one-dimensional shallow water equations in Lagrangian coordinates,
which admit symmetries and possess finite-difference analogues of the conservation laws of the original model.
On their basis, new schemes were also constructed for various extended models, such as the Green--Naghdi
equations~\cite{bk:DorKapMelGN2020} and the modified shallow water equations~\cite{bk:KaptDorMel_ModSW2022}. In the present paper, based on the results of~\cite{DorKap_2020,bk:DorKapSW_JMP_2020} and the paper~\cite{bk:ChevDorKap2020} devoted to the construction of schemes for wave equations, the authors construct and implement new finite-difference schemes for the one-dimensional SMHD equations in Lagrangian coordinates and in mass Lagrangian coordinates.

\medskip

The paper is organized as follows.
In Section~\ref{sec:01_eqns}, the SMHD equations in Eulerian and Lagrangian coordinates are considered and their dimensionless forms are given.
Sections~\ref{sec:02_schemes_lagr} and~\ref{sec:03_schemes_mass} are devoted to the construction of schemes for the SMHD equations in Lagrangian and mass Lagrangian coordinates and their conservation laws are listed for various bottom topographies according to the group classification results.
In Section~\ref{sec:04_invariance} invariance of finite-difference schemes in Lagrangian coordinates and invariance the constructed schemes in particular are discussed. Section~\ref{sec:05_numerical} is devoted to the numerical implementation of the schemes. The results are summarized in~Conclusion.

\section{The one-dimensional SMHD equations in Lagrangian coordinates}
\label{sec:01_eqns}

The SMHD model proposed in~\cite{bk:SWMHD_Gilman2000} in Eulerian coordinates has the form
\begin{subequations}
\label{GilmanSys}
\begin{equation}
h_{t}+\nabla^{\prime}\cdot(h\mathbf{u})=0,
\end{equation}
\begin{equation}
\mathbf{u}_{t}=-\nabla\left(\frac{\mathbf{u}\cdot\mathbf{u}}{2}\right)+\nabla\left(\frac{\mathbf{H}\cdot\mathbf{H}}{2}\right)-(\hat{\mathbf{k}}\times\mathbf{u})\hat{\mathbf{k}}\cdot\nabla\times\mathbf{u}+(\hat{\mathbf{k}}\times\mathbf{H})\hat{\mathbf{k}}\cdot\nabla\times\mathbf{H}-g\nabla h,
\end{equation}
\begin{equation}
\mathbf{H}_{t}=\nabla\times(\mathbf{u}\times\mathbf{H})+(\nabla^{\prime}\cdot\mathbf{u})\mathbf{H}-(\nabla^{\prime}\cdot\mathbf{H})\mathbf{u},
\end{equation}
\begin{equation}
\nabla^{\prime}\cdot(h\mathbf{H})=0,
\end{equation}
\end{subequations}
where $\mathbf{u}=(u,v,0)$ and $\mathbf{H}=(H^{x},H^{y},0)$
are the two-dimensional velocity and magnetic filed vectors,
$\hat{\mathbf{k}}=(0,0,1)$ denotes the unit vector in the vertical
direction, $\nabla^{\prime}\cdot$ is the horizontal divergence
operator, $\hat{\mathbf{k}}\cdot\nabla\times$ is the vertical component
of the curl operator, and the constant $g\neq0$ characterizes the
gravitational acceleration. Here $h = h_0 + \eta$ is the depth of the fluid, 
where $\eta$ characterizes a deviation of the free surface from the undisturbed level $h_0 > 0$. It is also considered~$|\eta| < h_0$. 

\medskip

By analogy with the standard shallow water equations, one introduces the function~$b(x, y)$
characterizing topography of the bottom~\cite{bk:Stoker1948,bk:KapMelDor_PhA2022}.
Then, one writes system~(\ref{GilmanSys}) with uneven bottom in coordinate form as
\begin{subequations}
\begin{equation}
h_{t}+uh_{x}+hu_{x}+vh_{y}+hv_{y}=0,
\end{equation}
\begin{equation}
u_{t}+uu_{x}+vu_{y}-H^{x}H_{x}^{x}-H^{y}H_{y}^{x}+gh_{x}=b_{x},
\end{equation}
\begin{equation}
v_{t}+uv_{x}+vv_{y}-H^{x}H_{x}^{y}-H^{y}H_{y}^{y}+gh_{y}=b_{y},
\end{equation}
\begin{equation}
H_{t}^{x}+uH_{x}^{x}+vH_{y}^{x}-u_{x}H^{x}-u_{y}H^{y}=0,
\end{equation}
\begin{equation}
H_{t}^{y}+uH_{x}^{y}+vH_{y}^{y}-v_{x}H^{x}-v_{y}H^{y}=0,
\end{equation}
\begin{equation}
h_{x}H^{x}+hH_{x}^{x}+h_{y}H^{y}+hH_{y}^{y}=0.
\end{equation}
\end{subequations}
Assuming that $h$, $b$, $\mathbf{u}$ and $\mathbf{H}$ only depend on the single space variable~$x$,
the latter system is brought to the form
\begin{subequations}
\label{eq:jul30.1}
\begin{equation}
h_{t}+uh_{x}+hu_{x}=0,\label{Eul_1d_cont}
\end{equation}
\begin{equation}
u_{t}+uu_{x}-H^{x}H_{x}^{x}+gh_{x}=b^{\prime},
\end{equation}
\begin{equation}
v_{t}+uv_{x}-H^{x}H_{x}^{y}=0,
\end{equation}
\begin{equation}
H_{t}^{x}+uH_{x}^{x}-u_{x}H^{x}=0,\label{Eul_1d_Hx}
\end{equation}
\begin{equation}
H_{t}^{y}+uH_{x}^{y}-v_{x}H^{x}=0,
\end{equation}
\begin{equation}
h_{x}H^{x}+hH_{x}^{x}=0.\label{Eul_divH}
\end{equation}
\end{subequations}
Using~(\ref{Eul_1d_cont}) and~(\ref{Eul_divH})
one can rewrite~(\ref{Eul_1d_Hx}) as
\[
0 = h_{t}H^{x}+hH_{t}^{x}=(hH^{x})_{t}.
\]
Thus, by means of (\ref{Eul_divH}),
\[
hH^{x}=a,
\]
where $a$ is constant that characterizes magnitude of the magnetic force.

\medskip

Further we consider the model in mass Lagrangian coordinates that is
\begin{subequations} \label{GilmanMassLagr}
\begin{equation}
\left(\frac{1}{h}\right)_t - u_s = 0,
\end{equation}
\begin{equation}
u_t - a^2 \left(\frac{1}{h}\right)_s - g h h_s - b^\prime = 0,
\end{equation}
\begin{equation}
h H^x = a,
\end{equation}
\begin{equation} \label{GilmanMassLagrInt1}
v_t - a H^y_s = 0,
\end{equation}
\begin{equation} \label{GilmanMassLagrInt2}
H^y_t - a v_s = 0,
\end{equation}
\end{subequations}
where the Lagrangian coordinates are introduced similarly to the gas dynamics equations \cite{bk:Ovsyannikov[2003]}
through the relations
\begin{equation}
\varphi_{t}(s,t)={u}(s,t),
\qquad
\varphi_{s}(s,t)=\frac{1}{{h}(s,t)}\label{eq:jul30.2}
\end{equation}
in such a way that the following relation
for the differentials~$dt$, $ds$ and~$dx$ holds~\cite{bk:YanenkRojd[1968]}
\[
ds=h\,dx-hu\,dt,
\]
and,
\[
s_{t}=-hu,\qquad s_{x}=h.
\]
The d'Alembert solution of the acoustic equations~(\ref{GilmanMassLagrInt1}), (\ref{GilmanMassLagrInt2}) has the form
\begin{equation}
v=f_{1}(s+at)+f_{2}(s-at),\qquad H^{y}=f_{1}(s+at)-f_{2}(s-at),\label{eq:aug03.1}
\end{equation}
where $f_{1}$ and $f_{2}$ are arbitrary functions of their arguments.
Thus, as they are integrated, we can exclude them from consideration when performing numerical modeling.

Notice that in the variables $\varphi$, $t$ and $s$ the system is reduced to the single equation~\cite{bk:KapMelDor_PhA2022}
\begin{equation} \label{GilmanLagr}
\varphi_{tt} - a^2 \varphi_{ss} + \left(\frac{g}{2 \varphi_s^2}\right)_s - b^\prime = 0,
\end{equation}
which is more suitable for constructing numerical schemes~\cite{DorKap_2020,bk:DorKapSW_JMP_2020,bk:DorKapMelGN2020,bk:KaptDorMel_ModSW2022}.

\medskip

The dimensionless form of system (\ref{GilmanMassLagr}) is the following
\begin{subequations}
\begin{equation}
\left(\frac{1}{\tilde{h}}\right)_{\tilde{t}} - \tilde{u}_{\tilde{s}} = 0,
\end{equation}
\begin{equation} \label{EvolUDimls}
\tilde{u}_{\tilde{t}}
    - \alpha^2 \left(\frac{1}{\tilde{h}}\right)_{\tilde{s}}
    - g_1 \tilde{h} \tilde{h}_{\tilde{s}} - \tilde{b}^\prime = 0,
\end{equation}
\begin{equation}
\tilde{v}_{\tilde{t}} - \alpha^2 \tilde{H}^y_{\tilde{s}} = 0,
\end{equation}
\begin{equation}
\tilde{H}^y_{\tilde{t}} - \tilde{v}_{\tilde{s}} = 0,
\end{equation}
\begin{equation}
\tilde{h} \tilde{H}^x = 1, \label{hHDimls}
\end{equation}
\end{subequations}
where $\alpha$ and $g_1$ are dimensionless constants characterizing intensity of magnetic and gravity fields.
Then, the solution (\ref{eq:aug03.1}) becomes
\[
\tilde{v}=\alpha f_{1}(\tilde{s}+\alpha \tilde{t})+\alpha f_{2}(\tilde{s}- \alpha \tilde{t}),
\qquad
\tilde{H}^{y}=f_{1}(\tilde{s}+\alpha \tilde{t})-f_{2}(\tilde{s}-\alpha \tilde{t}).
\]

\smallskip

The dimensionless form of (\ref{GilmanLagr}) is
\begin{equation} \label{GilmanLagrDimls}
\tilde{\varphi}_{\tilde{t}\tilde{t}}
    - \alpha^2 \tilde{\varphi}_{\tilde{s}\tilde{s}}
    + \left(\frac{g_1}{2 \tilde{\varphi}_{\tilde{s}}^2}\right)_{\tilde{s}}
    - {\tilde{b}}^\prime = 0.
\end{equation}
Further we consider the dimensionless equations, and symbol $\tilde{ }$ is omitted for brevity.
We also assume $g_1=2$ by means of equivalence transformations~\cite{bk:KapMelDor_PhA2022}.

\bigskip

The local conservation laws of~(\ref{GilmanLagrDimls}) have been obtained in~\cite{bk:KapMelDor_PhA2022}. They are listed below depending on the bottom topography according to the results of the group classifications with respect to the function~$b^\prime$.
\begin{itemize}
\item Case $b^\prime$ is arbitrary. Conservation laws of momentum and energy are
\[
\left(\varphi_{t}\varphi_{s}
\right)_{t}-\left(\frac{\varphi_{t}^{2}+{\alpha^2}\varphi_{s}^{2}}{2}-\frac{2}{\varphi_{s}}+b
\right)_{s}=0,
\]
\[
\left(\frac{\varphi_{t}^{2}+{\alpha^2}\varphi_{s}^{2}}{2}+\frac{1}{\varphi_{s}}
-b\right)_{t}+\left(
\frac{\varphi_{t}}{\varphi_{s}^{2}}-{\alpha^2}\varphi_{t}\varphi_{s}
\right)_{s}=0,
\]
and the conservation law of mass is just the relation $\varphi_{ts} = \varphi_{st}$;

\item Case $b^\prime=0$ (a horizontal bottom). Center-of-mass motion law
\[
\left(t\varphi_{t}-\varphi\right)_{t}+\left(\frac{t g}{2\varphi_{s}^{2}}-t\alpha^{2}\varphi_{s}\right)_{s}=0,
\]
and the following alternative form of the conservation law of momentum.
\[
\left(\varphi_{t}\right)_{t}+\left(\frac{1}{\varphi_{s}^{2}}-{\alpha^2}\varphi_{s}\right)_{s}=0;
\]

\item
The case of inclined bottom topography ($b^\prime=C$) is reduced to
the case of horizontal bottom topography by means of the
transformation~\cite{dorodnitsyn2019shallow,bk:SWMHD_Karelsky2014}
\begin{equation} \label{InclToFlat}
\varphi \mapsto \varphi + \frac{C t^2}{2};
\end{equation}

\item Case $b^\prime=\varphi$ (a parabolic bottom, concave up).
\[
\left(\left(\varphi-\varphi_{t}-\varphi_{s}\right)e^{t}\right)_{t}+\left(\left({\alpha^2}\varphi_{s}+\varphi_{t}+\varphi-\frac{1}{\varphi_{s}^{2}}\right)e^{t}\right)_{s}=0,
\]
\[
\left(\left(\varphi_{t}+\varphi_{s}+\varphi\right)e^{-t}\right)_{t}+\left(\left(\varphi-{\alpha^2}\varphi_{s}-\varphi_{t}+\frac{1}{\varphi_{s}^{2}}\right)e^{-t}\right)_{s}=0;
\]

\item Case $b^\prime=-\varphi$ (a parabolic bottom, concave down).
\[
\left(\varphi\sin{t}+\varphi_{t}\cos{t}\right)_{t}-\left(\left({\alpha^2}\varphi_{s}-\frac{1}{\varphi_{s}^{2}}\right)\cos{t}\right)_{s}=0,
\]
\[
\left(\varphi\cos{t}-\varphi_{t}\sin{t}\right)_{t}+\left(\left({\alpha^2}\varphi_{s}-\frac{1}{\varphi_{s}^{2}}\right)\sin{t}\right)_{s}=0.
\]
\end{itemize}

\section{Schemes for the SMHD equations in Lagrangian coordinates $\varphi(t,s)$}
\label{sec:02_schemes_lagr}

Further, to avoid confusion with the standard finite-difference notation, we denote $h=\rho$ and $\varphi = x$. 
The standard notation~\cite{bk:SamarskyPopov_book[1992]} is used where $\hat{f}$ and $\check{f}$
denote the values of~$f=f(t,s)$ at the point shifted up and down along the time axis. Similarly $f_+$ and $f_-$ (or $f^+$ and $f^-$) denote right and left shifts along the space axis.
The total differentiations are defined as follows
\[
f_t = \frac{1}{\tau} (\hat{f} - f),
\qquad
f_{\check{t}} = \frac{1}{\tau} (f - \check{f}),
\qquad
f_s = \frac{1}{h} (f_+ - f),
\qquad
f_{\bar{s}} = \frac{1}{h} (f - f_-).
\]
This notation should not cause confusion with standard partial derivatives, since the rest of the discussion is dedicated to finite differences.

\medskip

The authors would like to emphasize here that the main goal of the present publication is the construction of symmetry-preserving finite-difference schemes that possess the largest possible number of conservation laws. The study of such issues as convergence and stability of the constructed schemes, their well-posedness and regularity of solutions goes beyond the scope of our work. The interested reader can find a detailed discussion of nonlinear stability of schemes for hyperbolic systems and related issues, e.g., in~\cite{bk:NonlinearStabilityBook}.

\medskip

Considering the form of~(\ref{GilmanLagrDimls}), one notices that in case $\alpha=0$ (i.e., magnetic field is absent) one gets the one-dimensional shallow water equations in Lagrangian coordinates (see, e.g., \cite{bk:KaptsovMeleshko_1D_classf[2018]}).
Finite-difference schemes for the shallow water equations have been constructed by the authors in their recent publications~\cite{DorKap_2020,bk:DorKapSW_JMP_2020}.
Here, equation~(\ref{GilmanLagrDimls}) differs from the shallow water equations only by a linear term $\alpha^2\varphi_{ss}$, which can be easily approximated.
The construction can be accomplished by combining the schemes previously obtained in~\cite{DorKap_2020}
and the scheme for the linear wave equation derived in~\cite{bk:ChevDorKap2020}.
Thus, one writes the following finite-difference scheme on an orthogonal uniform mesh
\begin{subequations} \label{SchemeMain}
\begin{equation}
\displaystyle
x_{t\check{t}} - {\alpha}^2 x_{s\bar{s}} + \left(\frac{1}{\hat{x}_s \check{x}_s}\right)_{\bar{s}} - \check{B} = 0,
\end{equation}
\begin{equation}
h_+ = h_- = h,
\qquad
\hat{\tau} = \check{\tau} = \tau,
\qquad
\vec{h} \cdot \vec{\tau} = 0,
\end{equation}
\end{subequations}
where $h$ and $\tau$ are steps along the space and time axes, and $\check{B}$ is some approximation for $b^\prime$.
Scheme~(\ref{SchemeMain}) has the second order of approximation in~$h$ and~$\tau$.

\bigskip

The conservation laws for the scheme are derived by straightforward extension of the previous results~\cite{DorKap_2020,
bk:ChevDorKap2020}. Recall that any (local) conservation law of (\ref{SchemeMain}) can be represented in the divergent form
\begin{equation} \label{CLLambdaForm}
 T_t + S_s = \Lambda \left(
    x_{t\check{t}} - {\alpha}^2 x_{s\bar{s}} + \left(\frac{1}{\hat{x}_s \check{x}_s}\right)_{\bar{s}} - \check{B}
 \right) = 0,
\end{equation}
where $T$ and $S$ are conserved quantities of the conservation law, and $\Lambda$ is called a conservation law multiplier.

\bigskip

The finite-difference counterparts of the conservation laws listed in Section~\ref{sec:02_schemes_lagr} possessed by scheme~(\ref{SchemeMain}) and the corresponding conservation law multipliers are listed below. 
One can verify by direct calculations that these conservation laws are satisfied on the solutions of scheme~(\ref{SchemeMain}) (and therefore can also be represented in form~(\ref{CLLambdaForm})). To do this, it suffices, for example, to solve the equations of the scheme with respect to $x_t$, $h_+$, and $\hat{\tau}$ and substitute the result and its finite-difference shifts into the conservation laws.
\begin{itemize}

    \item
    In case $b(x)=\text{const}$ ($\check{B}=0$), the conservation laws are the following.
    \begin{enumerate}
     \item
     mass
        \[
            (\hat{x}_s)_{\check{t}} - (x_t^+)_{\bar{s}} = 0;
        \]

      \item
      momentum
      \[
        \Lambda_1 = 1, \qquad
        (x_t)_{\check{t}} + \left(
            (\hat{x}_s \check{x}_s)^{-1}
            - {\alpha}^2 x_s
        \right)_{\bar{s}} = 0;
      \]

      \item
      center-of-mass law
      \[
        \Lambda_2 = t, \qquad
        (t x_t - x)_{\check{t}} + \left(
            t (\hat{x}_s \check{x}_s)^{-1}
            - t {\alpha}^2 x_s
        \right)_{\bar{s}} = 0;
      \]

      \item
      energy
      \[
        \Lambda_3 = \frac{x_t + \check{x}_t}{2}, \qquad
        \frac{1}{2} (
            x_t^2 + x_s^{-1} + \hat{x}_s^{-1}
            + {\alpha}^2 x_s \hat{x}_s
        )_{\check{t}}
        + \frac{1}{2} \left(
            (x_t^+ + \check{x}_t^+) \left(
            (\hat{x}_s \check{x}_s)^{-1}
            - {\alpha}^2 x_s
            \right)
        \right)_{\bar{s}} = 0.
      \]
    \end{enumerate}

\item
The case of inclined bottom ($\check{B}=C$) is reduced to the case
of horizontal bottom topography by means of the following
finite-difference analogue of~(\ref{InclToFlat})
(see~\cite{dorodnitsyn2019shallow}).
\begin{equation} \label{InclToHrz}
x \mapsto x + \frac{C t \hat{t}}{2}.
\end{equation}

\item
For the parabolic bottom $\displaystyle b(x) = \frac{x^2}{2}$ (see also~\cite{bk:DorKapSW_JMP_2020}),
\[
\check{B} = \frac{2 (\cosh \tau - 1)}{\tau^2} x
\]
there are two additional conservation laws
\[
\Lambda_{{a}1}^+ = e^t,
\qquad
\left(
    x \frac{e^{\hat{t}} - e^t}{\tau}
    - e^t {x}_t
\right)_{\check{t}}
    -\left(
        e^t ((\hat{x}_s \check{x}_s)^{-1} - {\alpha}^2 x_s)
    \right)_{\bar{s}} = 0,
\]
\[
\Lambda_{{a}2}^+ = e^{-t},
\qquad
\left(
    x \frac{e^{-t} - e^{-\hat{t}}}{\tau}
    - e^{-t} {x}_t
\right)_{\check{t}}
    - \left(
        e^{-t} ((\hat{x}_s \check{x}_s)^{-1} - {\alpha}^2 x_s)
    \right)_{\bar{s}}
    = 0.
\]

\item
For the parabolic bottom $\displaystyle b(x) = -\frac{x^2}{2}$ (see also~\cite{bk:DorKapSW_JMP_2020}),
\[
\check{B} = \frac{2 (\cos \tau - 1)}{\tau^2} x
\]
there are two additional conservation laws
\[
\Lambda_{{a}1}^- = \sin t,
\qquad
\left(
    {x}_t \sin t
    - x \frac{\sin \hat{t} - \sin t}{\tau}
\right)_{\check{t}}
    + \left(
        \sin t \, ((\hat{x}_s \check{x}_s)^{-1} - {\alpha}^2 x_s)
    \right)_{\bar{s}}
 = 0,
\]
\[
\Lambda_{{a}2}^- = \cos t,
\qquad
\left(
    {x}_t \cos t
    - x \frac{\cos \hat{t} - \cos t}{\tau}
\right)_{\check{t}}
    + \left(
        \cos t \, ((\hat{x}_s \check{x}_s)^{-1} - {\alpha}^2 x_s)
    \right)_{\bar{s}}
    = 0.
\]

\item
In case the bottom is arbitrary ($b=b(x)$), one can preserve the conservation law of energy by modifying the
first equation of the scheme as follows (see also~\cite{DorKap_2020}).
\begin{equation} \label{scmArbBtn}
x_{t\check{t}} - {\alpha}^2 x_{s\bar{s}} + \left(\frac{1}{\hat{x}_s \check{x}_s}\right)_{\bar{s}}
    - \frac{b_t + \check{b}_t}{x_t + x_{\check{t}}}  = 0.
\end{equation}
Then, the conservation law of energy becomes
 \[
        (
            x_t^2 + x_s^{-1} + \hat{x}_s^{-1}
            + {\alpha}^2 x_s \hat{x}_s
            - \hat{b} - b
        )_{\check{t}}
        + \left(
            (x_t^+ + \check{x}_t^+) \left(
            (\hat{x}_s \check{x}_s)^{-1}
            - {\alpha}^2 x_s
            \right)
        \right)_{\bar{s}} = 0.
\]

\end{itemize}

\section{Schemes for the SMHD equations in mass Lagrangian coordinates}
\label{sec:03_schemes_mass}

Instead of the three-layer scheme in Lagrangian coordinates one can consider a two-layer scheme in mass Lagrangian coordinates. This can be achieved by introducing a specific approximation of the `state equation' $p = \rho^2$
and for the transformation~(\ref{eq:jul30.2}) as it was previously done in~\cite{DorKap_2020,bk:DorKapSW_JMP_2020,bk:DorKapMelGN2020,bk:KaptDorMel_ModSW2022}.
The resulting scheme is
\begin{equation} \label{mass_lagr_sq}
    \def\arraystretch{2}
    \begin{array}{c}
    \displaystyle
    \rho_{\check{t}}
        + \frac{1}{2} \rho \check{\rho} \left(u_s + \check{u}_s\right) = 0,
    \\
    \displaystyle
      u_{\check{t}} + Q_{\bar{s}} - \check{B} = 0,
    \\
    \displaystyle
    x_t = u, \qquad
    \check{x}_s + x_s = \frac{1}{\sqrt{\check{p}}} + \frac{1}{\sqrt{p}} = \frac{2}{\check{\rho}},
    \\
    h^s_+ = h^s_- = h, \quad
    \hat{\tau} = \check{\tau} = \tau,
    \quad
    \vec{h} \cdot \vec{\tau} = 0,
    \end{array}
\end{equation}
where the quantity $Q$ is given by
\begin{equation}\label{flux_Q}
  Q = \left(\frac{4}{\rho\check{\rho}}
    - \frac{2}{\sqrt{p}}\left( \frac{1}{\rho} + \frac{1}{\check{\rho}} \right) + \frac{1}{p}
    \right)^{-1}
    - \frac{{\alpha}^2}{\sqrt{p}},
\end{equation}
and the equation
\begin{equation} \label{mass_lagr_sq_state}
\frac{1}{\sqrt{\check{p}}} + \frac{1}{\sqrt{p}} = \frac{2}{\check{\rho}}
\end{equation}
approximates the equation~$p = \rho^2$.
Scheme (\ref{mass_lagr_sq}) along with equation (\ref{mass_lagr_sq_state}) have the first order of approximation in~$h$ and~$\tau$.

\medskip
In mass Lagrangian coordinates, the conservation laws are the following.
\begin{itemize}
\item
In case $b(x)=\text{const}$ ($\check{B}=0$):
\begin{enumerate}
    \item
    mass
    \[
        \left( \frac{1}{\rho} \right)_{\check{t}}  - \left(
            \frac{u^+ + \check{u}^+}{2}
        \right)_{\bar{s}} = 0;
    \]
    \item
    momentum
    \[
    u_{\check{t}}  + Q_{\bar{s}} = 0;
    \]
    \item
    center-of-mass law
    \[
    (t u - x)_{\check{t}}  + \left(t Q \right)_{\bar{s}} = 0;
    \]
    \item
    energy
    \[
        \displaystyle
        \left(
            \frac{u^2}{2} + \frac{p}{2 \sqrt{p} - \rho}
            + \frac{{\alpha}^2 (2 \sqrt{p} - \rho)}{2 p \rho}
        \right)_{\check{t}}
        + \left(
            \frac{u^+ + \check{u}^+}{2} \, Q
        \right)_{\bar{s}} = 0.
    \]

    \end{enumerate}

    \item
    As mentioned above, the case of an inclined bottom is reduced to the case of a horizontal bottom.

    \item
    For the parabolic bottom $\displaystyle b(x) = \frac{x^2}{2}$:
    \[
        \left(
        x \frac{e^{\hat{t}} - e^t}{\tau}
        - e^t u
        \right)_{\check{t}}
        -\left(
        e^t Q
        \right)_{\bar{s}}
        = 0,
    \]
    \[
        \left(
        x \frac{e^{-t} - e^{-\hat{t}}}{\tau}
        - e^{-t} u
        \right)_{\check{t}}
        - \left(
        e^{-t} Q
        \right)_{\bar{s}}
        = 0.
    \]

    \item
    For the parabolic bottom $\displaystyle b(x) = -\frac{x^2}{2}$:
    \[
    \left(
        x \frac{
            \sin \hat{t} - \sin t}{\tau}
            - u \sin t
        \right)_{\check{t}}
        - \left(Q \sin t\right)_{\bar{s}} = 0,
    \]
    \[
        \left(
            x \frac{\cos \hat{t} - \cos t}{\tau}
            - u \cos t
        \right)_{\check{t}}
        - \left(
            Q \cos t
        \right)_{\bar{s}} = 0.
    \]

\item
In case the bottom is arbitrary ($b=b(x)$), the conservation law of energy is
\begin{equation} \label{clEnergyInMassCoordsB}
        \displaystyle
        \left(
            \frac{u^2}{2} + \frac{p}{2 \sqrt{p} - \rho}
            + \frac{{\alpha}^2 (2 \sqrt{p} - \rho)}{2 p \rho}
            - \frac{\hat{b} + b}{2}
        \right)_{\check{t}}
        + \left(
            \frac{u^+ + \check{u}^+}{2} \, Q
        \right)_{\bar{s}} = 0.
    \end{equation}

\end{itemize}
As in the case of scheme~(\ref{SchemeMain}), the validity of the listed conservation laws on solutions of~(\ref{mass_lagr_sq}) can be verified by direct calculations.

\section{Discussion on the invariance of the constructed schemes}
\label{sec:04_invariance}

In the previous sections we focused on the conservation laws of the constructed
finite-difference schemes. Here we mention that the constructed schemes are not
only conservative~(i.e., possessing conservation laws), but also invariant,
that is, they preserve the symmetries of the original differential model.

Recall that there is a close relationship between conservation laws and symmetries of an equation, which can be expressed in the form of Noether's theorem~\cite{bk:Noether1918,bk:Ibragimov2[2007]}.
In the finite-difference case for ordinary difference equations, there is a finite-difference analogue of Noether's theorem~\cite{bk:DorodnitsynKozlovWint[2003], bk:Dorodnitsyn[2011]}, which makes it possible to derive conservation laws from known symmetries admitted by the equation. In the case of finite-difference schemes for partial differential equations, conservation laws are usually established using the direct
method~\cite{DorKap_2020,bk:ChevDorKap2020} or algebraic transformations~\cite{bk:SamarskyPopov_book[1992], bk:SamarskyPopov[1970], bk:DorKapMDPI2022,bk:KapDorMel_MHDCylSchemes2023,bk:KaptsovMHDCylSchemes_Conf2022}. In the latter case, knowledge of the symmetries and conservation laws of the original differential model may suggest the form of finite-difference conservation laws~\cite{bk:DorKapMDPI2022,bk:KapDorMel_MHDCylSchemes2023,bk:KaptsovMHDCylSchemes_Conf2022}.
As far as the authors know, there is no rigorous proof of the converse statement, i.e., from the existence of a conservation law for a scheme, in general, it does not follow that the scheme admits a symmetry corresponding to the conservation law. At the same time, a large number of examples~\cite{DorKap_2020,bk:DorKapSW_JMP_2020,bk:ChevDorKap2020,bk:DorKapMelGN2020,bk:KaptDorMel_ModSW2022,bk:DorKapMDPI2022,bk:KapDorMel_MHDCylSchemes2023} indirectly indicate the existence of such a connection. 
In particular, as can be verified, scheme~(\ref{SchemeMain}) for various bottom topographies admits the same sets of symmetries (Lie algebras) as equation~(\ref{GilmanLagrDimls}). Here we refer to~\cite{bk:KapMelDor_PhA2022} where the Lie algebras admitted by~(\ref{GilmanLagrDimls}) are stated as a result of the group classification procedure.

\medskip

The second significant factor that should be taken into account when constructing invariant finite-difference schemes is the preservation of the orthogonality and uniformity of the finite-difference mesh by group transformations. The corresponding mesh invariance criteria are given in~\cite{Dorodnitsyn1991,bk:Dorodnitsyn[2011]}.
In contrast to Eulerian coordinates, in Lagrangian coordinates in the one-dimensional case, these criteria are usually satisfied (this was the case for the Lie algebras admitted by the schemes for shallow water equations~\cite{DorKap_2020}, modified shallow water equations~\cite{bk:KaptDorMel_ModSW2022}, Green--Naghdi equations~\cite{bk:DorKapMelGN2020}, and this is also true for the obtained SMHD schemes~(\ref{SchemeMain}) and~(\ref{mass_lagr_sq})). Thus, the advantage of using Lagrangian coordinates for one-dimensional equations is the possibility of constructing symmetry-preserving schemes on simple orthogonal uniform meshes. We note that in higher-dimensional cases the situation becomes more complicated and some symmetries, such as relabeling symmetries, may not satisfy the criteria as it was stated in~\cite{bk:DorKapMel_SW2D}. 
As simplicial meshes are often used in practical applications for numerical calculations in two and three spatial dimensions, the orthogonality requirement in this case turns out to be quite severe. It can be weakened with the loss of some symmetries (for example, relabeling symmetries seem not to be generally admitted by simplicial meshes). Although such meshes may naturally arise in the construction of schemes admitting certain Lie algebras,  
this is primarily refers to the finite element or finite volume method, the discussion of which is beyond the scope of our study.

\section{Numerical implementation of the constructed schemes}
\label{sec:05_numerical}

In the present section, some typical problems are solved numerically using schemes of the form~(\ref{mass_lagr_sq}) in mass Lagrangian coordinates.
The numerical implementation is carried out using pseudo-viscosity, which often leads to a loss of accuracy at discontinuities, but which is quite enough to evaluate the qualitative picture of solutions. Scheme~(\ref{mass_lagr_sq}) is linearized using the Newton method~(see, e.g., \cite{bk:SamarskyPopov_book[1992]}), after which it is transformed to a form suitable for solving through an iterative algorithm and the tridiagonal matrix method~\cite{bk:Samarskii2001theory}. Since the linearized version of the scheme may include derivatives of approximations of the function describing the bottom topography, we first discuss the possible forms of these approximations.

\subsection{Remark on approximations for the function describing the bottom topography}

As it was mentioned above, $\check{B}$ is some approximation of the derivative~$b^\prime$ of the function~$b$ that describes the bottom topography. For the linearization purposes, we need to know the set of variables on which approximations~$\check{B}$ may depend. According to the group classification~\cite{bk:KapMelDor_PhA2022}, there the following possibilities are of interest: $b=C$, $b=C x$, $b=\pm\frac{x^2}{2}$, $b=\ln x$, and $b=b(x)$. As in~(\ref{scmArbBtn}), we consider approximations
\[
\check{B} = \frac{[b(x) + b(\check{x})]_t}{x_t + \check{x}_t}.
\]
The approximations for the listed topographies are the following.
\begin{enumerate}[1)]
  \item
  In case $b=C$ one gets $\check{B}=0$.

  \item
  In case $b=C x$ one gets $\check{B} = C$.

  \item
  In case $b=\pm\frac{x^2}{2}$,
  \[
  \check{B}
    = \pm{x} \pm \frac{\tau^2}{2} \check{u}_t
    = \pm\frac{x + \check{x}}{2} \pm \frac{\tau}{2} u.
  \]

  \item
  In case $b=\ln x$,
  \[
  \check{B} =  \frac{1}{\tau (u + \check{u})} \, \ln \frac{x + \tau u}{x - \tau \check{u}}
    = \frac{1}{x + \tau u - \check{x}} \ln \frac{x + \tau u}{\check{x}}
    = \frac{1}{x} + O(\tau^2).
  \]
\end{enumerate}
Thus, in general we can restrict our consideration to the function $\check{B} = \check{B}(x, \check{x}, u)$, or
\[
B = B(\hat{x}, x, \hat{u}).
\]

\subsection{Linearization of the schemes}
\label{sec:linearization}

We utilize Newton's linearization method, for which we first augment scheme~(\ref{mass_lagr_sq}) 
with the equation~\cite{bk:SamarskyPopov_book[1992]}
\[
\omega - \Omega(\rho, u_s) = 0,
\]
where 
\[
    \Omega(\rho, u_s) = \frac{1}{2} \nu \rho (u_s - |u_s|)
    + \frac{1}{2} \mu \rho |u_s| (u_s - |u_s|).
\]
The grid function $\omega$ introduces linear and quadratic pseudo-viscosities characterized by the corresponding coefficients~$\nu$ and~$\mu$.
When implementing scheme~(\ref{mass_lagr_sq}), for the equation of motion, instead of~(\ref{flux_Q}) we also consider the quantity~$\Theta = Q - \omega.$

Linearizing the scheme, one derives the equations
\begin{equation*}
\def\arraystretch{2}
\begin{array}{r}
    \displaystyle
    \left(
        1 + \frac{\tau}{2 h} \left(\itr{k}{u}_{i+1} - \itr{k}{u}_i + u^j_{i+1} - u^j_i\right) \rho^j_i
    \right) \left(
        \itr{k+1}{\rho}_i - \itr{k}{\rho}_i
    \right)
    + \frac{\tau}{2 h} \rho^j_i \itr{k}{\rho}_i
        \left( \itr{k+1}{u_{i+1}} - \itr{k}{u_{i+1}} - \itr{k+1}{u_{i}} + \itr{k}{u_{i}} \right)
    = -\itr{k}{f_{1,i}},
    \\
        \displaystyle
        \left(1 - \tau \frac{\partial}{\partial{\itr{k}{u_i}}}{B\left( \itr{k}{x_i}, x^j_i, \itr{k}{u_i} \right)} \right)
        \left(\itr{k+1}{u_{i}} - \itr{k}{u_{i}}\right)
        - \tau \frac{\partial}{\partial{\itr{k}{x_i}}}{B\left(\itr{k}{x_i}, x^j_i, \itr{k}{u_i} \right)}
        \left(\itr{k+1}{x_{i}} - \itr{k}{x_{i}}\right)
        \\
        \displaystyle
        + \frac{\tau}{h}
            \left( \itr{k+1}{\Theta_{i}} - \itr{k}{\Theta_{i}} - \itr{k+1}{\Theta_{i-1}} + \itr{k}{\Theta_{i-1}} \right)
        = -\itr{k}{f_{2,i}},
    \\
    \displaystyle
    \itr{k+1}{x_{i}} - \itr{k}{x_{i}}  = -\itr{k}{f_{3,i}},
    \\
    \displaystyle
    \itr{k+1}{P_{i}} - \itr{k}{P_{i}}  = -\itr{k}{f_{4,i}},
\end{array}
\end{equation*}
\begin{equation*}
\begin{array}{r}
    \begin{cases}
        \displaystyle
        \itr{k+1}{\omega_i} - \itr{k}{\omega_i}
        + \frac{1}{h} \left(
            \nu \left( \itr{k}{u_{i+1}} - \itr{k}{u_{i}} \right)
            - \mu \left( \itr{k}{u_{i+1}} - \itr{k}{u_{i}} \right)^2
        \right) (\itr{k+1}{\rho_{i}} - \itr{k}{\rho_{i}})   & \\
        \displaystyle
        \qquad + \frac{1}{h} \itr{k}{\rho_{i}} \left(
            \nu - 2 \mu \left( \itr{k}{u_{i+1}} - \itr{k}{u_{i}} \right)
        \right) \left(
            \itr{k+1}{u_{i+1}}
            - \itr{k}{u_{i+1}}
            - \itr{k+1}{u_{i}}
            + \itr{k}{u_{i}}
        \right)
        = -\itr{k}{f_{5,i}}, & \mbox{if } \itr{k}{u_{i+1}} < \itr{k}{u_{i}} \\
      \itr{k+1}{\omega_i} - \itr{k}{\omega_i} = -\itr{k}{f_{5,i}}, & \mbox{otherwise};
    \end{cases}
\end{array}
\end{equation*}
\begin{equation*}
\begin{array}{r}
    \displaystyle
    \itr{k+1}{\Theta_i} - \itr{k}{\Theta_i}
    - \itr{k+1}{\omega_i} + \itr{k}{\omega_i}
    + \left(\itr{k+1}{\rho_{i}} - \itr{k}{\rho_{i}}\right)
    \left(
        \frac{2 \left(\itr{k}{P_i}\right)^3 \rho^j_i}
        {(\rho^j_i - 2 \itr{k}{P_i})(\itr{k}{\rho_i} - 2 \itr{k}{P_i})^2}
    \right)
    \\
    \displaystyle
    \quad
    + \left(\itr{k+1}{P_i} - \itr{k}{P_i}\right)
    \left(
        \frac{2 \itr{k}{P_i} \rho^j_i \itr{k}{\rho_i} \left(
                (\rho^j_i + \itr{k}{\rho_i})\itr{k}{P_i}
                - \rho^j_i \itr{k}{\rho_i}
            \right)
        }
        {(\rho^j_i - 2 \itr{k}{P_i})^2(\itr{k}{\rho_i} - 2 \itr{k}{P_i})^2}
        - \frac{\alpha^2}{\left(\itr{k}{P_i}\right)^2}
    \right)
    = -\itr{k}{f_{6,i}},
\end{array}
\end{equation*}
where $P = p^2$, $\itr{k}{q}$ denotes the value of the quantity $q$ on the $k$-th iteration, and
\begin{equation*}
\begin{array}{c}
\displaystyle
\itr{k}{f_{1,i}}
    = \itr{k}{\rho_i} - \rho^j_i + \frac{\tau}{2 h} \rho^j_i \itr{k}{\rho_i} \left(
        \itr{k}{u_{i+1}} - \itr{k}{u_{i}}
        + u^j_{i+1} - u^j_i
    \right),
\\
\displaystyle
\itr{k}{f_{2,i}}
    = \itr{k}{u_i} - u^j_i + \frac{\tau}{h} \left( \itr{k}{\Theta_i} - \itr{k}{\Theta_{i-1}} \right)
        - \tau B\left( \itr{k}{x_i}, x^j_i, \itr{k}{u_i} \right),
\\
\displaystyle
\itr{k}{f_{3,i}}
    = \itr{k}{x_i} - x^j_i - \tau u^j_i,
\\
\displaystyle
\itr{k}{f_{4,i}}
    = \itr{k}{P_i} - \frac{\rho^j_i P^j_i}{2 P^j_i - \rho^j_i},
\\
\displaystyle
\itr{k}{f_{5,i}}
    =
    \begin{cases}
        \displaystyle
      \itr{k}{\omega_i}
        + \frac{\nu}{h} \itr{k}{\rho_i} \left(\itr{k}{u_{i+1}} - \itr{k}{u_{i}}\right)
         -\frac{\mu}{h} \itr{k}{\rho_i} \left(\itr{k}{u_{i+1}} - \itr{k}{u_{i}}\right)^2, & \mbox{if } \itr{k}{u_{i+1}} < \itr{k}{u_{i}} \\
      \itr{k}{\omega_i}, & \mbox{otherwise};
    \end{cases}
\\
\displaystyle
\itr{k}{f_{6,i}}
    = \itr{k}{\Theta_i} - \itr{k}{\omega_i}
    - \left[
        \frac{4}{\rho^j_i \itr{k}{\rho_i}}
        - \frac{2}{\itr{k}{P_i}}\left(
            \frac{1}{\itr{k}{\rho_i}}
            + \frac{1}{\rho^j_i}
        \right)
        + \frac{1}{\left(\itr{k}{P_i}\right)^2}
    \right]^{-1}
    + \frac{\alpha^2}{\itr{k}{P_i}}.
\end{array}
\end{equation*}

Notice that in case $\mu^2 + \nu^2 \neq 0$ conservation law~(\ref{clEnergyInMassCoordsB}) cannot be preserved.
Therefore, when verifying the conservation law of energy, one should take the results of calculations obtained without pseudo-viscosity or with very small values of $\nu$ and $\mu$.

\medskip

Introducing variables
\[
\itr{k+1}{\delta u} = \itr{k+1}{u} - \itr{k}{u},
\quad
\itr{k+1}{\delta \rho} = \itr{k+1}{\rho} - \itr{k}{\rho},
\quad
\itr{k+1}{\delta x} = \itr{k+1}{x} - \itr{k}{x},
\quad
\cdots,
\]
one derives the three-point equation
\[
A_i \itr{k+1}{\delta u_{i-1}}
- C_i \itr{k+1}{\delta u_{i}}
+ D_i \itr{k+1}{\delta u_{i+1}}
= - F_i,
\]
where $A_i$, $B_i$, $D_i$, and $F_i$ depend on the values obtained on the previous iterations.
The latter equation can be solved using tridiagonal matrix algorithm~\cite{bk:Samarskii2001theory,bk:SamarskyPopov_book[1992]}.
Then, provided $\itr{k+1}{\delta u}$, one can calculate the values of~$\itr{k+1}{\delta \rho}$, $\itr{k+1}{\delta x}$, etc. In particular,
\[
\displaystyle
\itr{k+1}{\delta \rho_i} = -\frac{\tau \itr{k}{\rho_i} \check{\rho_i} \left(\delta u_{i+1} - \delta u_i\right) + 2 h \itr{k}{f_{1,i}} }
        {\tau \check{\rho}_i \left(\itr{k}{u}_{i+1} - \itr{k}{u}_i + \check{u}_{i+1} - \check{u}_i\right) + 2 h}.
\]

\medskip

Further we consider only problems in which velocity of the flow at the boundaries is held constant. Such problems correspond to the simplest left (L) and right (R) boundary conditions of the form~\cite{bk:SamarskyPopov_book[1992]} $\delta{u}_L=\delta{u}_R=0$.

\subsection{Numerical calculations}

First, we consider the dam break problem over horizontal, parabolic and logarithmic bottoms. At the initial moment of time, the liquid is divided in the middle of the computational domain $0 \leqslant s \leqslant S$ into two regions with heights $\rho_L$ and $\rho_R < \rho_L$, and the liquid height profile is slightly smoothed near the initial discontinuity. We chose the following calculation parameters.
\[
\tau = 0.05 h,
\quad
S = 4.0,
\quad
\rho_L = 1.0,
\quad
\rho_R = 0.5,
\quad
\alpha^2 = 1.6,
\]
\[
\nu = \nu_0 h,
\quad
\mu = \frac{3}{2\pi^2} \mu_0^2 h^2.
\]
The viscosity coefficients $\nu_0 = 1 \div 2$, $\mu_0 = 3 \div 4$,
and the step $h = 0.04 \div 0.1$ are vary depending on the problem and the bottom topography under consideration.

By means of equations (\ref{EvolUDimls}) and (\ref{hHDimls}), one can approximately write
\[
\hat{u} \sim u + \tau \alpha^2 \frac{\partial{H^x}}{\partial{s}} + \cdots.
\]
This means that the change in the velocity of the flow depends linearly on the slope of the profile of the longitudinal component $H^x$ of the magnetic field.
Further, we are only interested in the shape of this profile, therefore, to visualize the profile, we introduce the quantity
\begin{equation} \label{gradHx}
B = \kappa \left(\frac{1}{h}\right)_s = \kappa H^x_s,
\end{equation}
where the coefficient $\kappa$ is selected for each problem so that $B$ changes within close to the area of change in the liquid height~$\rho$.

\begin{figure}[H]
  \centering
  \includegraphics[width=0.75\linewidth]{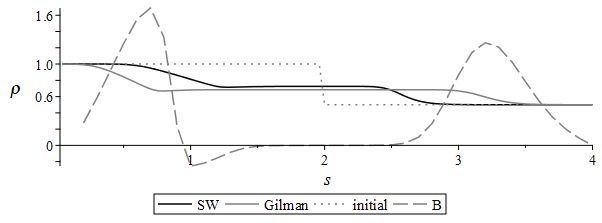}
  \caption{Dam break problem for the horizontal bottom. `SW' (black) is the solution for the standard shallow water equations ($\alpha=0$), `Gilman' is the solution for the case $\alpha^2 > 0$ (SMHD). Initial profile is given by the dotted line, and the magnetic field gradient $B$ is denoted by the dashed line.}
  \label{fig:fig1}
\end{figure}

The obtained solutions are depicted in Figure~\ref{fig:fig1}, Figure~\ref{fig:fig2} and Figure~\ref{fig:fig3} for the time $t = 0.92$.
The parabolic bottom profile is given by the function
\[
b(x) = 0.5 k_p (x - 0.5 S)^2, \qquad k_p = 0.05,
\]
and the logarithmic bottom profile is described by
\[
b(x) = k_{l}^1 \ln(x + k_{l}^2),
\qquad
k_{l}^1 = 0.1,
\quad
k_{l}^2 = 2.0.
\]

\begin{figure}[H]
  \centering
  \includegraphics[width=0.5\linewidth]{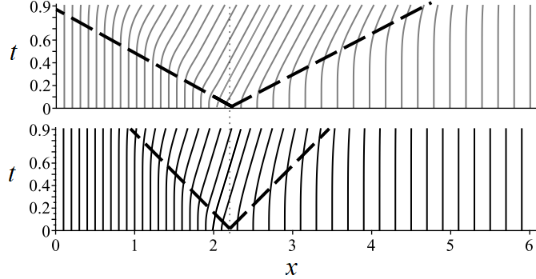}
  \caption{Dam break problem for the horizontal bottom. Trajectories $x(t,s)$ for the case $\alpha > 0$ (gray) and $\alpha = 0$ (black). 
    The characteristics of the flow are outlined with thick dashed lines.}
  \label{fig:fig1traj}
\end{figure}

\begin{figure}[H]
  \centering
  \includegraphics[width=0.75\linewidth]{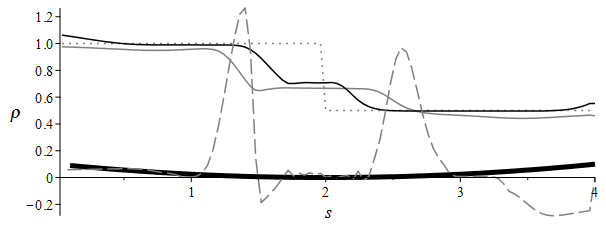}
  \caption{Dam break problem for the parabolic bottom.}
  \label{fig:fig2}
\end{figure}

\begin{figure}[H]
  \centering
  \includegraphics[width=0.75\linewidth]{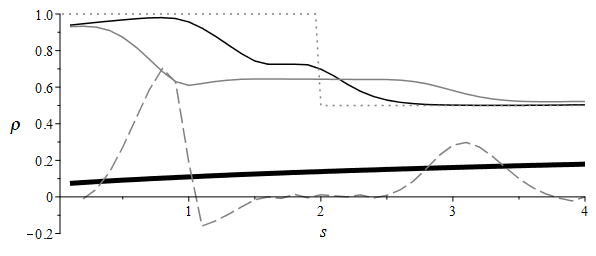}
  \caption{Dam break problem for the logarithmic bottom.}
  \label{fig:fig3}
\end{figure}
\begin{figure}[H]
  \centering
  \includegraphics[width=0.75\linewidth]{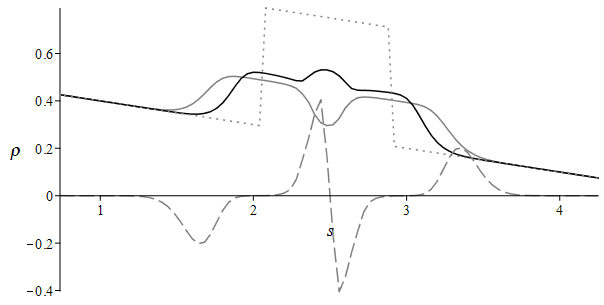}
  \caption{Collapse of the liquid column over the inclined bottom.}
  \label{fig:fig4}
\end{figure}

The second problem is related to the calculation of the collapse of a liquid column over the inclined bottom (Figure~\ref{fig:fig4}) given by the function
\[
b(x) = k_{\text{incl}} \, x, \qquad
k_{\text{incl}} = -0.1.
\]
In Figure~\ref{fig:fig3}, the initial column profile is shown as a dotted gray line.
One sees that in the presence of a magnetic field the column collapses faster (gray solid line) than in the case of standard shallow water model (black solid line). Notice that the same result is obtained by calculating for a horizontal bottom followed by applying transformation~(\ref{InclToHrz}).

\smallskip

In all the experiments, an acceleration of fluid motion in the areas of growth of the gradient of~$H^x$ is observed (see also trajectories of motion depicted in~Figure~\ref{fig:fig1traj}), as predicted by formula~(\ref{gradHx}). 
In case $b^\prime =0$, the shock wave velocity
\[
\mathcal{D} = \sqrt{(\alpha/h_0)^2 + g_1 h_0}
\]
is obtained using the Rankine--Hugoniot conditions (see, e.g., \cite{Bouchut2016} for details). Here $h_0$ is the height of the liquid column in front of the wave. Accordingly, the characteristics of the flow are given in Figure~\ref{fig:fig1traj}.

Notice that the difference analogue (\ref{mass_lagr_sq}) of~(\ref{eq:jul30.2}),
\[
\displaystyle
\frac{\check{x}_s + x_s}{2} - \frac{1}{2\sqrt{\check{p}}} - \frac{1}{2\sqrt{p}} = 0,
\]
which has the series expansion $\displaystyle x_s(t,s) - \frac{1}{\rho(t, s)} = O(h + \tau) $, 
remains correct, since according to Figures~\ref{fig:fig1} and \ref{fig:fig2}---\ref{fig:fig4} 
the density~$\rho$ satisfies the inequality $0 < \rho < \infty$ in the computational domain.

Finally, 
we estimate the change in the total energy with time for the first problem using the following error estimations
\[
err_{A}(n) = |\mathcal{H}(n) - \mathcal{H}(0)|,
\qquad
err_{R}(n) = \frac{|\mathcal{H}(n) - \mathcal{H}(0)|}{|\mathcal{H}(0)|},
\]
where, according to the conservation law~(\ref{clEnergyInMassCoordsB}), 
\[
\displaystyle
\mathcal{H}(n) = h \sum_{k=0}^{\lfloor S/h \rfloor} \left(
            \frac{(u^n_k)^2}{2} + \frac{p^n_k}{2 \sqrt{p^n_k} - \rho^n_k}
            + \frac{{\alpha}^2 (2 \sqrt{p^n_k} - \rho^n_k)}{2 p^n_k \rho^n_k}
            - b(x^n_k)
    \right),
\]
and the index notation is used: $f^n_k = f(t + n \tau, s + k h) $.
$\mathcal{H}(n)$ gives the total energy in the computational domain, and its value should
tend to constant in the continuous limit~\cite{bk:Dorod_Hamilt[2011]}.

\begin{figure}[H]
  \centering
  \includegraphics[width=0.75\linewidth]{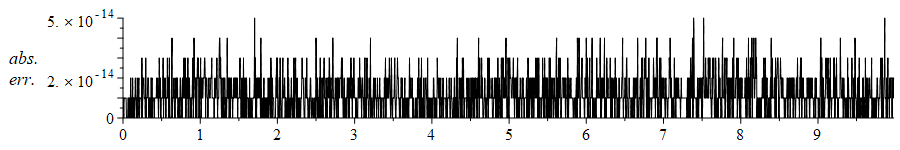}
  \vfill
  \includegraphics[width=0.75\linewidth]{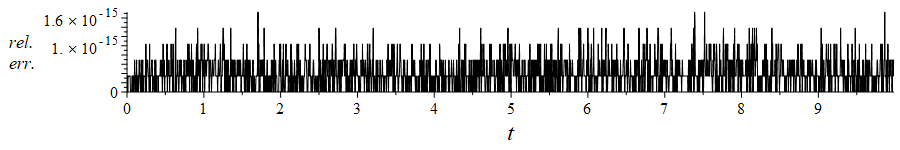}
  \caption{Total energy preservation estimations $err_A(n(t))$ and $err_R(n(t))$, $n(t) = \lfloor t / \tau \rfloor$ for $t \in [0,10]$. The calculation is performed without pseudo-viscosity. Machine precision is $10^{-16}$.}
  \label{fig:cl}
\end{figure}

The result for the time interval $0 \leqslant t \leqslant  10$ is given in Figure~\ref{fig:cl}.
The figure demonstrates that the total energy practically does not change with time.

\section{Conclusion}

In the present paper, finite-difference schemes have been constructed
for one-dimensional shallow water equations in the presence of a
magnetic field (the Gilman model) in Lagrangian coordinates for
various bottom topographies. According to the recent group
classification~\cite{bk:KapMelDor_PhA2022}, the cases of an
arbitrary, horizontal, inclined, parabolic, and logarithmic bottom
are distinguished, while the case of an inclined bottom is reduced
to the case of a horizontal bottom by a simple point change of
variables. The constructed schemes possess finite-difference
analogues of the conservation laws of the original differential
model for all the listed cases of bottom topography, and are also
invariant, i.e. preserve the symmetries of the original model. The
schemes are constructed on uniform orthogonal meshes. In Lagrangian
coordinates, they are given on three time layers, and in mass
Lagrangian coordinates, they can be given on two time layers by
means of a specially selected finite-difference equation of state.

In mass Lagrangian coordinates, the schemes are implemented
numerically. Typical one-dimensional problems of dam break and
liquid column collapse above an inclined bottom are considered. To
demonstrate the effect of a magnetic field on the motion of fluid
particles, the Gilman model is compared to the standard shallow
water model. Numerical experiments show that the magnetic field
accelerates the movement of compression waves along the various
bottom topographies, and the destruction of the liquid column above
the inclined bottom. For one of the problems calculated with no
pseudo-viscosity, the control of the conservation law of energy was
carried out, from which it is seen that the invariant conservative
scheme preserves the total energy almost with no loss.

\subsection*{Acknowledgements}

This work is a continuation of research supported by
the Russian Science Foundation Grant No. 18-11-00238 ``Hydrodynamics-type equations:
symmetries, conservation laws, invariant difference schemes''.
The authors are grateful to S. V. Meleshko for valuable discussions. 
E.I.K.~also sincerely appreciates the hospitality of the Suranaree University of Technology.


\label{lastpage}
\end{document}